\documentclass{amsart}
\usepackage{amsfonts}

\newtheorem{theorem}{Theorem}
\theoremstyle{plain}

\newtheorem{remark}{Remark}

\numberwithin{equation}{section}
\input{tcilatex}

\begin{document}
\title[Ostrowski Inequality]{A Generalisation of an Ostrowski Inequality in
Inner Product Spaces}
\author{Sever S. Dragomir}
\address{School of Computer Science and Mathematics\\
Victoria University of Technology\\
PO Box 14428, MCMC 8001\\
Victoria, Australia.}
\email{sever@matilda.vu.edu.au}
\urladdr{http://rgmia.vu.edu.au/SSDragomirWeb.html}
\author{Anca C. Go\c{s}a}
\address{College No. 12, Re\c{s}i\c{t}a, Jud. Cara\c{s}-Severin, R0-1700, Re%
\c{s}i\c{t}a, Romania.}
\email{ancagosa@hotmail.com}

\begin{abstract}
A generalisation of inner product spaces of an inequality due to Ostrowski
and applications for sequences and integrals are given.
\end{abstract}

\date{May 19, 2003}
\maketitle

\section{Introduction}

In 1951, A.M. Ostrowski \cite[p. 289]{2b} obtained the following result (see
also \cite[p. 92]{1b}).

\begin{theorem}
\label{t1.1}Suppose that $\mathbf{a}=\left( a_{1},\dots ,a_{n}\right) ,$ $%
\mathbf{b}=\left( b_{1},\dots ,b_{n}\right) $ and $\mathbf{x}=\left(
x_{1},\dots ,x_{n}\right) $ are real $n-$tuples such that $\mathbf{a}\neq 0$
and%
\begin{equation}
\sum_{i=1}^{n}a_{i}x_{i}=0\text{ \ and \ }\sum_{i=1}^{n}b_{i}x_{i}=1.
\label{1.1}
\end{equation}%
Then%
\begin{equation}
\sum_{i=1}^{n}x_{i}^{2}\geq \frac{\sum_{i=1}^{n}a_{i}^{2}}{%
\sum_{i=1}^{n}a_{i}^{2}\sum_{i=1}^{n}b_{i}^{2}-\left(
\sum_{i=1}^{n}a_{i}b_{i}\right) ^{2}},  \label{1.2}
\end{equation}%
with equality if and only if%
\begin{equation}
x_{k}=\frac{b_{k}\sum_{i=1}^{n}a_{i}^{2}-a_{k}\sum_{i=1}^{n}a_{i}b_{i}}{%
\sum_{i=1}^{n}a_{i}^{2}\sum_{i=1}^{n}b_{i}^{2}-\left(
\sum_{i=1}^{n}a_{i}b_{i}\right) ^{2}},\ \ k=1,\dots ,n.  \label{1.3}
\end{equation}
\end{theorem}

Another similar result due to Ostrowski which is far less known and obtained
in the same work \cite[p. 130]{2b} (see also \cite[p. 94]{1b}), is the
following one.

\begin{theorem}
\label{t1.2}Let $\mathbf{a}$, $\mathbf{b}$ and $\mathbf{x}$ be $n-$tuples of
real numbers with $\mathbf{a}\neq 0$ and%
\begin{equation}
\sum_{i=1}^{n}a_{i}x_{i}=0\text{ \ and \ }\sum_{i=1}^{n}x_{i}^{2}=1.
\label{1.4}
\end{equation}%
Then%
\begin{equation}
\frac{\sum_{i=1}^{n}a_{i}^{2}\sum_{i=1}^{n}b_{i}^{2}-\left(
\sum_{i=1}^{n}a_{i}b_{i}\right) ^{2}}{\sum_{i=1}^{n}a_{i}^{2}}\geq \left(
\sum_{i=1}^{n}b_{i}x_{i}\right) ^{2}.  \label{1.5}
\end{equation}%
If $\mathbf{a}$ and $\mathbf{b}$ are not proportional, then the equality
holds in (\ref{1.5}) iff%
\begin{equation}
x_{k}=q\cdot \frac{b_{k}\sum_{i=1}^{n}a_{i}^{2}-a_{k}\sum_{i=1}^{n}a_{i}b_{i}%
}{\left( \sum_{k=1}^{n}a_{k}^{2}\right) ^{\frac{1}{2}}\left[
\sum_{i=1}^{n}a_{i}^{2}\sum_{i=1}^{n}b_{i}^{2}-\left(
\sum_{i=1}^{n}a_{i}b_{i}\right) ^{2}\right] ^{\frac{1}{2}}},\ \ k\in \left\{
1,\dots ,n\right\} ,  \label{1.6}
\end{equation}%
with $q\in \left\{ -1,1,\right\} .$
\end{theorem}

The case of equality which was neither mentioned in \cite{1b} nor in \cite%
{2b} is considered in Remark \ref{r1}.

In the present paper, by the use of an elementary argument based on\
Schwarz's inequality, a natural generalisation in inner-product spaces of (%
\ref{1.5}) is given. The case of equality is analyzed. Applications for
sequences and integrals are also provided.

\section{The Results}

The following theorem holds.

\begin{theorem}
\label{t2.1}Let $\left( H,\left\langle \cdot ,\cdot \right\rangle \right) $
be a real or complex inner product space and $a,b\in H$ two linearly
independent vectors. If $x\in H$ is such that%
\begin{equation}
\left\langle x,a\right\rangle =0\ \text{and\ }\left\Vert x\right\Vert =1, 
\tag{(i)}  \label{(i)}
\end{equation}%
then%
\begin{equation}
\frac{\left\Vert a\right\Vert ^{2}\left\Vert b\right\Vert ^{2}-\left\vert
\left\langle a,b\right\rangle \right\vert ^{2}}{\left\Vert a\right\Vert ^{2}}%
\geq \left\vert \left\langle x,b\right\rangle \right\vert ^{2}.  \label{2.1}
\end{equation}%
The equality holds in (\ref{2.1}) iff%
\begin{equation}
x=\nu \left( b-\frac{\overline{\left\langle a,b\right\rangle }}{\left\Vert
a\right\Vert ^{2}}\cdot a\right) ,  \label{2.2}
\end{equation}%
where $\nu \in \mathbb{K}$ $\left( \mathbb{C},\mathbb{R}\right) $ is such
that%
\begin{equation}
\left\vert \nu \right\vert =\frac{\left\Vert a\right\Vert }{\left[
\left\Vert a\right\Vert ^{2}\left\Vert b\right\Vert ^{2}-\left\vert
\left\langle a,b\right\rangle \right\vert ^{2}\right] ^{\frac{1}{2}}}.
\label{2.3}
\end{equation}
\end{theorem}

\begin{proof}
We use Schwarz's inequality in the inner product space $H,$ i.e.,%
\begin{equation}
\left\Vert u\right\Vert ^{2}\left\Vert v\right\Vert ^{2}\geq \left\vert
\left\langle u,v\right\rangle \right\vert ^{2},\ \ u,v\in H  \label{2.4}
\end{equation}%
with equality iff there is a scalar $\alpha \in \mathbb{K}$ such that%
\begin{equation}
u=\alpha v.  \label{2.5}
\end{equation}%
If we apply (\ref{2.4}) for $u=z-\frac{\left\langle z,c\right\rangle }{%
\left\Vert c\right\Vert ^{2}}\cdot c,$ $v=d-\frac{\left\langle
d,c\right\rangle }{\left\Vert c\right\Vert ^{2}}\cdot c$, where $c\neq 0$
and $c,d,z\in H,$ and taking into account that%
\begin{align*}
\left\Vert z-\frac{\left\langle z,c\right\rangle }{\left\Vert c\right\Vert
^{2}}\cdot c\right\Vert ^{2}& =\frac{\left\Vert z\right\Vert ^{2}\left\Vert
c\right\Vert ^{2}-\left\vert \left\langle z,c\right\rangle \right\vert ^{2}}{%
\left\Vert c\right\Vert ^{2}}, \\
\left\Vert d-\frac{\left\langle d,c\right\rangle }{\left\Vert c\right\Vert
^{2}}\cdot c\right\Vert ^{2}& =\frac{\left\Vert d\right\Vert ^{2}\left\Vert
c\right\Vert ^{2}-\left\vert \left\langle d,c\right\rangle \right\vert ^{2}}{%
\left\Vert c\right\Vert ^{2}}
\end{align*}%
and%
\begin{equation*}
\left\langle z-\frac{\left\langle z,c\right\rangle }{\left\Vert c\right\Vert
^{2}}\cdot c,d-\frac{\left\langle d,c\right\rangle }{\left\Vert c\right\Vert
^{2}}\cdot c\right\rangle =\frac{\left\langle z,d\right\rangle \left\Vert
c\right\Vert ^{2}-\left\langle z,c\right\rangle \left\langle
c,d\right\rangle }{\left\Vert c\right\Vert ^{2}},
\end{equation*}%
we deduce the inequality%
\begin{equation}
\left[ \left\Vert z\right\Vert ^{2}\left\Vert c\right\Vert ^{2}-\left\vert
\left\langle z,c\right\rangle \right\vert ^{2}\right] \left[ \left\Vert
d\right\Vert ^{2}\left\Vert c\right\Vert ^{2}-\left\vert \left\langle
d,c\right\rangle \right\vert ^{2}\right] \geq \left\vert \left\langle
z,d\right\rangle \left\Vert c\right\Vert ^{2}-\left\langle z,c\right\rangle
\left\langle c,d\right\rangle \right\vert ^{2}  \label{2.6}
\end{equation}%
with equality iff there is a $\beta \in \mathbb{K}$ such that%
\begin{equation}
z=\frac{\left\langle z,c\right\rangle }{\left\Vert c\right\Vert ^{2}}\cdot
c+\beta \left( d-\frac{\left\langle d,c\right\rangle }{\left\Vert
c\right\Vert ^{2}}\cdot c\right) .  \label{2.7}
\end{equation}%
If in (\ref{2.6}) we choose $z=x,$ $c=a$ and $d=b,$ where $a$ and $x$
statisfy (i), then we deduce%
\begin{equation*}
\left\Vert a\right\Vert ^{2}\left[ \left\Vert a\right\Vert ^{2}\left\Vert
b\right\Vert ^{2}-\left\vert \left\langle a,b\right\rangle \right\vert ^{2}%
\right] \geq \left[ \left\langle x,b\right\rangle \left\Vert a\right\Vert
^{2}\right] ^{2}
\end{equation*}%
which is clearly equivalent to (\ref{2.1}).

The equality holds in (\ref{2.1}) iff%
\begin{equation*}
x=\nu \left( b-\frac{\overline{\left\langle a,b\right\rangle }}{\left\Vert
a\right\Vert ^{2}}\cdot a\right) ,
\end{equation*}%
where $\nu \in \mathbb{K}$ satisfies the condition%
\begin{equation}
1=\left\Vert x\right\Vert =\left\vert \nu \right\vert \left\Vert b-\frac{%
\overline{\left\langle a,b\right\rangle }}{\left\Vert a\right\Vert ^{2}}%
\cdot a\right\Vert =\left\vert \nu \right\vert \left[ \frac{\left\Vert
a\right\Vert ^{2}\left\Vert b\right\Vert ^{2}-\left\vert \left\langle
a,b\right\rangle \right\vert ^{2}}{\left\Vert a\right\Vert ^{2}}\right] ^{%
\frac{1}{2}},  \label{2.8}
\end{equation}%
and the theorem is thus proved.
\end{proof}

The following particular cases hold.

\begin{enumerate}
\item[\textbf{1.}] If $\mathbf{a}$, $\mathbf{b}$, $\mathbf{x}\in \ell
^{2}\left( \mathbb{K}\right) ,$ $\mathbb{K}=\mathbb{C},\mathbb{R},$ where 
\begin{equation*}
\ell ^{2}\left( \mathbb{K}\right) :=\left\{ x=\left( x_{i}\right) _{i\in 
\mathbb{N}},\ \sum_{i=1}^{\infty }\left\vert x_{i}\right\vert ^{2}<\infty
\right\}
\end{equation*}%
with $\mathbf{a}$, $\mathbf{b}$ linearly independent and%
\begin{equation}
\sum_{i=1}^{\infty }x_{i}\overline{a_{i}}=0,\ \ \ \sum_{i=1}^{\infty
}\left\vert x_{i}\right\vert ^{2}=1,  \tag{a}  \label{a}
\end{equation}%
then%
\begin{equation}
\frac{\sum_{i=1}^{\infty }\left\vert a_{i}\right\vert ^{2}\sum_{i=1}^{\infty
}\left\vert b_{i}\right\vert ^{2}-\left\vert \sum_{i=1}^{\infty }a_{i}%
\overline{b_{i}}\right\vert ^{2}}{\sum_{i=1}^{\infty }\left\vert
a_{i}\right\vert ^{2}}\geq \left\vert \sum_{i=1}^{\infty }x_{i}\overline{%
b_{i}}\right\vert ^{2}.  \label{2.9}
\end{equation}%
The equality holds in (\ref{2.9}) iff%
\begin{equation}
x_{i}=\nu \left[ b_{i}-\frac{\sum_{k=1}^{\infty }a_{k}\overline{b_{k}}}{%
\sum_{k=1}^{\infty }\left\vert a_{k}\right\vert ^{2}}\cdot a_{i}\right] ,\ \
\ i\in \left\{ 1,2,\dots \right\}  \label{2.10}
\end{equation}%
with $\nu \in \mathbb{K}$ is such that%
\begin{equation}
\left\vert \nu \right\vert =\frac{\left( \sum_{k=1}^{\infty }\left\vert
a_{k}\right\vert ^{2}\right) ^{\frac{1}{2}}}{\left[ \sum_{k=1}^{\infty
}\left\vert a_{k}\right\vert ^{2}\sum_{k=1}^{\infty }\left\vert
b_{k}\right\vert ^{2}-\left\vert \sum_{k=1}^{\infty }a_{k}\overline{b_{k}}%
\right\vert ^{2}\right] ^{\frac{1}{2}}}.  \label{2.11}
\end{equation}
\end{enumerate}

\begin{remark}
\label{r1}The case of equality in (\ref{1.5}) is obviously a particular case
of the above. We omit the details.
\end{remark}

\begin{enumerate}
\item[\textbf{2.}] If $f,g,h\in L^{2}\left( \Omega ,m\right) ,$ where $%
\Omega $ is an $m-$measurable space and 
\begin{equation*}
L^{2}\left( \Omega ,m\right) :=\left\{ f:\Omega \rightarrow \mathbb{K},\
\int_{\Omega }\left\vert f\left( x\right) \right\vert ^{2}dm\left( x\right)
<\infty \right\} ,
\end{equation*}%
with $f,g$ being linearly independent and%
\begin{equation}
\int_{\Omega }h\left( x\right) \overline{f\left( x\right) }dm\left( x\right)
=0,\ \ \ \int_{\Omega }\left\vert h\left( x\right) \right\vert ^{2}dm\left(
x\right) =1,  \label{2.12}
\end{equation}%
then%
\begin{multline}
\frac{\int_{\Omega }\left\vert f\left( x\right) \right\vert ^{2}dm\left(
x\right) \int_{\Omega }\left\vert g\left( x\right) \right\vert ^{2}dm\left(
x\right) -\left\vert \int_{\Omega }f\left( x\right) \overline{g\left(
x\right) }dm\left( x\right) \right\vert ^{2}}{\int_{\Omega }\left\vert
f\left( x\right) \right\vert ^{2}dm\left( x\right) }  \label{2.13} \\
\geq \left\vert \int_{\Omega }h\left( x\right) \overline{g\left( x\right) }%
dm\left( x\right) \right\vert ^{2}.
\end{multline}%
The equality holds in (\ref{2.13}) iff%
\begin{equation*}
h\left( x\right) =\nu \left[ g\left( x\right) -\frac{\int_{\Omega }g\left(
x\right) \overline{f\left( x\right) }dm\left( x\right) }{\int_{\Omega
}\left\vert f\left( x\right) \right\vert ^{2}dm\left( x\right) }f\left(
x\right) \right] \text{ \ for \ a.e. \ }x\in \Omega
\end{equation*}%
and\ $\nu \in \mathbb{K}$ with%
\begin{equation*}
\left\vert \nu \right\vert =\frac{\left( \int_{\Omega }\left\vert f\left(
x\right) \right\vert ^{2}dm\left( x\right) \right) ^{\frac{1}{2}}}{\left[
\int_{\Omega }\left\vert f\left( x\right) \right\vert ^{2}dm\left( x\right)
\int_{\Omega }\left\vert g\left( x\right) \right\vert ^{2}dm\left( x\right)
-\left\vert \int_{\Omega }f\left( x\right) \overline{g\left( x\right) }%
dm\left( x\right) \right\vert ^{2}\right] ^{\frac{1}{2}}}.
\end{equation*}
\end{enumerate}


\begin{thebibliography}{9}
\bibitem{1b} D.S. MITRINOVI\'{C}, J.E. PE\v{C}ARI\'{C} and\ A.M. FINK, 
\textit{Classical and New Inequalities in Analysis}, Kluwer Academic
Publishers, Dordrecht/Boston/London, 1993.

\bibitem{2b} A.M. OSTROWSKI, \textit{Varlesungen \"{u}ber Differential und
Integralrechnung II}, Birkh\"{a}user, Basel, 1957.
\end{thebibliography}
\end{document}